# LOGARITHMIC BARRIER OPTIMIZATION PROBLEM USING NEURAL NETWORK

A. K. Ojha, C. Mallick, D. Mallick

**Abstract**— The combinatorial optimization problem is one of the important applications in neural network computation. The solutions of linearly constrained continuous optimization problems are difficult with an exact algorithm, but the algorithm for the solution of such problems is derived by using logarithm barrier function. In this paper we have made an attempt to solve the linear constrained optimization problem by using general logarithm barrier function to get an approximate solution. In this case the barrier parameters behave as temperature decreasing to zero from sufficiently large positive number satisfying convexity of the barrier function. We have developed an algorithm to generate decreasing sequence of solution converging to zero limit.

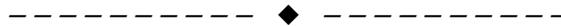

## 1 INTRODUCTION

Let us consider graph G= (V,E) where V=(1,2,...n) is the set of nodes and E is the set of edges. The edge between two nodes and is defined by (i,j). Let $(W_{ij})_{n \times n}$ is symmetric weight matrix such that

$$W_{ij} > 0 \ \forall \in E \qquad (1.1)$$

Given two nodes s and t, if s<t there exist a partition U and V-U to the node set V such that

$$W(s) = \sum_{i \in U} w_{ij} \qquad (1.2)$$

Where $j \in V - U$ is maximized for s, t not belonging to the same set. This problem is called s-t max cut problem. Through the exact solution is difficult however the same can be obtained by considering a equivalent problem.

$$Max: \tfrac{1}{4}\sum_{i=1}^{n}\sum_{j=1}^{n}(1 - x_i x_j)w_{ij} \qquad (1.3)$$

Subject to $x_s + x_t = 0$ and $x_i \in \{-1,1\}, \forall i = 1,2,\ldots n$

Benson and Zhang [2] found an algorithm to get an approximate solution in this direction to solve combinatorial and quadratic optimization problem.

Neural computation is combinatorial and quadratic optimization problem due to Hopfield and Tank [10] put a foundation stone in solving such type of problems. Many computational model and combinatorial optimization algorithm have been developed. Aiyer, Niranjan and Fallside [1] extended this idea for a theoretical study of such problem. Similar algorithm is used in solving travelling salesman problem by Durbin and Willshaw [6] and multiple traveling salesman problem by using Neural Network, due to Wacholder,

————————————————
- A.K. Ojha is with the School of Basic Sciences, Indian Institute of Technology, Bhubaneswar, Orissa-751014, India.
- C. Mallick is with the Department of Mathematicsics, B.O.S.E, Cuttack, Orissa-753007, India.
- D. Mallick is with the Department of Mathematics, Centurion Institute of Technology, Bhubaneswar, Orissa-751050, India.

© 2009 Journal of Computing

Han and Mann [15].Gee, Aiyer and Prager [7] made a frame work for studying optimization problem by Neural Network, where in a different paper Gee and Prager [7] used polyhedral combinatorics in Neural Network. Urahama [13] used gradient projection Network for solving linear and nonlinear programming problems. Some of the related works can be found in Bout and Miller[3],Gold and Rangaranjan [9],Gold and Mjolsness [11], Simic [12], Waugh and Westervelt [16],Walfe , Parry and Macmillan[17], Xu [18],Yuille and Kosowsky [19] etc. A systematic study of neural computational model is given by Van Den Berg [14], Clchoki and Unbehaunen [4] helps to formulate an algorithm to solve such type of problem in an efficient manner. All these work based on deterministic annealing type and to study the global minimum of the effective energy at high temperature and track it as temperature decreases. In this paper we have generalized the work of Dang [5] by considering an equivalent linear constrained optimizing problem with general logarithmic barrier function.

$$b(x) = -\sum_{i=1}^{n}[\ln(px_i + q) + \ln(q - px_i)] \quad (1.4)$$

to find an approximate solution. Where the parameters p, q are positive real. The barrier parameters behave as temperature decreasing to zero from a sufficiently large positive number satisfying convexity of the barrier function. We have developed an algorithm for which the network converges to at least a local minimum point, if a local minimum point of the barrier problem is generated for a decreasing sequence of values with zero limit. The organization of the paper is as follows. Following the introduction in section 2, we have considered logarithm barrier function to minimize the problem with Lagrange multiplier. In section 3, we have proved the theorem by considering general logarithmic function. Numerical example in section 4, has been presented to show the efficiency of the algorithm. Finally the paper concluded with a remark in the section 5.

## 2 Linear Logarithmic Barrier function:

In order to find the solution of the problem first of all we will convert s-t max cut problem into a linearly constrained continuous optimization problem and then derive some theoretical results from an application of linear logarithmic barrier function. From s-t max-cut problem we have

$$\min: f(x) = \frac{1}{2}x^T(W - \alpha I)x \quad (2.1)$$

Subject to $x_i + x_t = 0$

$x_i \in \{-1,1\}, i = 1,2,3,\ldots,n$

where $\alpha$ is any positive number and I an n×n identity matrix and

$$W = \begin{pmatrix} 0 & w_{12} & \cdots & w_{1n} \\ \vdots & \ddots & & \vdots \\ w_{n1} & w_{n2} & \cdots & 0 \end{pmatrix}$$

be a symmetric weight matrix.

Let $b(x) = -\sum_{i=1}^{n}[\ln(px_i + q) + \ln(q - px_i)]$ (2.2)

which is used as barrier function $-q/p < x_i < q/p$, i=1, 2,……n in the objective function.

Let us define a function $h(x, \beta)$ for a positive barrier parameter $\beta$ such that

$$\min h(x,\beta) = f(x) + \beta b(x) \quad (2.3)$$

Subject to $x_s + x_t = 0$

Let $\Gamma = \{x: x_s + x_t = 0\}$

and $B = \left\{x: -\frac{q}{p} \le x_i \le \frac{q}{p}\right\}, i = 1,2,\ldots n$

$$\frac{\partial b(x)}{\partial x_i} = -\left[\frac{p}{px_i + q} + \frac{-p}{q - px_i}\right] = -p\left[\frac{1}{px_i + q} - \frac{1}{q - px_i}\right]$$

Thus $\lim_{x_i \to (-\frac{q}{p})} \frac{\partial b(x)}{\partial x_i} = -\infty$, $\lim_{x_i \to (\frac{q}{p})} \frac{\partial b(x)}{\partial x_i} = \infty$ (2.4)

From the definition of f(x) it is clear that $\nabla f(x) = (W + \alpha I)x$ is bounded on a set B and hence there exists an interior point $x^*$ of B such that it is a point of local or global minimum point of $h(x, \beta)$ defined in (2.1).

Now we will develop an algorithm for approximating a solution of

Min: $f(x) = \frac{1}{2}x^T(W + \alpha I)x$

Subject to $x_s + x_t = 0$

$x_i \in \{-1,1\}, i = 1,2,\ldots,n$ and to find minimum point where every component equal to -1 and 1 from the solution of (2.1) at the limiting value of β which converging to zero.

Let $L(x,\lambda) = h(x,\beta) - \lambda e_{st}^T x$

$e_{st} = (0,0,\ldots,0,1,0,0,\ldots,0,1,0,\ldots,0)^T \in R^n$

For any given barrier parameter β>0, the first order necessary optimality condition of (2.1) says that if $x$ is a minimum point of (2.1) then there exists a Lagrange multiplier λ satisfying $\nabla_x L(x,\lambda) = \nabla h(x,\beta) - \lambda e_{st} = 0$

$e_{st}^T x = 0$ (2.5)

For $i \neq s$ or t we have

$\frac{\partial L(x,\lambda)}{\partial x_i} = \frac{\partial f(x)}{\partial x_i} - \beta \left(\frac{p}{px_i+q} - \frac{p}{q-px_i}\right) = 0$

$\Rightarrow \beta p \left[\frac{1}{px_i+q} - \frac{1}{q-px_i}\right] = \frac{\partial f(x)}{\partial x_i}$

$\Rightarrow \beta p \left[\frac{-2px_i}{q^2-(px_i)^2}\right] = \frac{\partial f(x)}{\partial x_i}$

$\Rightarrow (q)^2 \frac{\partial f(x)}{\partial x_i} - (px_i)^2 \frac{\partial f(x)}{\partial x_i} = -2\beta p \cdot px_i$

$\Rightarrow (px_i)^2 \frac{\partial f(x)}{\partial x_i} - 2\beta p \cdot px_i - (q)^2 \frac{\partial f(x)}{\partial x_i} = 0$

$\Rightarrow px_i = \frac{2\beta p \pm \sqrt{4\beta^2 p^2 - 4\frac{\partial f(x)}{\partial x_i}(-q^2\frac{\partial f(x)}{\partial x_i})}}{2\frac{\partial f(x)}{\partial x_i}}$

$\Rightarrow x_i = \frac{\beta p \pm \sqrt{\beta^2 p^2 + q^2\left(\frac{\partial f(x)}{\partial x_i}\right)^2}}{p\frac{\partial f(x)}{\partial x_i}}$ (2.6)

and for i=s or t

$\frac{\partial L(x,\lambda)}{\partial x_i} = \left(\frac{\partial f(x)}{\partial x_i} - \lambda\right) - \beta \left(\frac{p}{px_i+q} - \frac{p}{q-px_i}\right) = 0$

We derive $x_i = \frac{\beta p \pm \sqrt{\beta^2 p^2 + q^2\left(\frac{\partial f(x)}{\partial x_i} - \lambda\right)^2}}{p\left(\frac{\partial f(x)}{\partial x_i} - \lambda\right)}$

Let $d_i(x,\lambda) = \begin{cases} \frac{\beta p \pm \sqrt{\beta^2 p^2 + q^2\left(\frac{\partial f(x)}{\partial x_i}\right)^2}}{p\left(\frac{\partial f(x)}{\partial x_i}\right)}, & i \neq s \text{ or } t \\ \frac{\beta p \pm \sqrt{\beta^2 p^2 + q^2\left(\frac{\partial f(x)}{\partial x_i} - \lambda\right)^2}}{p\left(\frac{\partial f(x)}{\partial x_i} - \lambda\right)}, & i = s \text{ or } t \end{cases}$ (2.7)

For i=1, 2,…n and

$d(x,\lambda) = (d_1(x,\lambda), d_2(x,\lambda), \ldots d_n(x,\lambda))^T$ from the above discussion it shows that $(d_i(x,\lambda) - x$ is a descent direction of $L(x,\lambda)$ if $-\frac{q}{p} < x < \frac{q}{p}$. Before proving our main theorem we will prove the following lemma.

## 2.1. Lemma

Assume that $-\frac{q}{p} < x < \frac{q}{p}$

(i) If $d_i(x,\lambda) - x_i < 0$ then $\frac{\partial}{\partial x_i}(L(x,\lambda)) > 0$

(ii) If $d_i(x,\lambda) - x_i > 0$ then $\frac{\partial}{\partial x_i}(L(x,\lambda)) < 0$

(iii) If $d_i(x,\lambda) - x_i = 0$ then $\frac{\partial}{\partial x_i}(L(x,\lambda)) = 0$

(iv) If $d_i(x,\lambda) - x_i \neq 0$ then $\nabla_x L(x,\lambda)^T(d(x,\lambda) - x) < 0$

(v) If $e_{st}^T(d(x,\lambda) - x) = 0$ and $d(x,\lambda) - x \neq 0$ then $\nabla h(x,\beta)^T(d(x,\lambda) - x) < 0$

**Proof:** Before proving these lemma we need to prove that $\left(\frac{\partial L(x,\lambda)}{\partial x_i}\right) > 0$ if $d_i(x,\lambda) - x_i < 0$

Let $v_i = \begin{cases} 0, & i \neq s \text{ or } t \\ \lambda, & i = s \text{ or } t \end{cases}$

$i = 1,2,\ldots n$

When $\frac{\partial f(x)}{\partial x_i} - v_i \neq 0$

$\frac{\beta p \pm \sqrt{\beta^2 p^2 + q^2\left(\frac{\partial f(x)}{\partial x_i} - v_i\right)^2}}{p\left(\frac{\partial f(x)}{\partial x_i} - v_i\right)} =$

$\frac{\left\{\beta p + \sqrt{\beta^2 p^2 + q^2\left(\frac{\partial f(x)}{\partial x_i} - v_i\right)^2}\right\}\left\{\beta p - \sqrt{\beta^2 p^2 + q^2\left(\frac{\partial f(x)}{\partial x_i} - v_i\right)^2}\right\}}{p\left(\frac{\partial f(x)}{\partial x_i} - v_i\right)\left\{\beta p - \sqrt{\beta^2 p^2 + q^2\left(\frac{\partial f(x)}{\partial x_i} - v_i\right)^2}\right\}}$

$= \frac{\beta^2 p^2 - \beta^2 p^2 - q^2\left(\frac{\partial f(x)}{\partial x_i} - v_i\right)^2}{p\left(\frac{\partial f(x)}{\partial x_i} - v_i\right)\left\{\beta p - \sqrt{\beta^2 p^2 + q^2\left(\frac{\partial f(x)}{\partial x_i} - v_i\right)^2}\right\}}$

$= \frac{-q^2\left(\frac{\partial f(x)}{\partial x_i} - v_i\right)^2}{p\left(\frac{\partial f(x)}{\partial x_i} - v_i\right)\left\{\beta p - \sqrt{\beta^2 p^2 + q^2\left(\frac{\partial f(x)}{\partial x_i} - v_i\right)^2}\right\}}$

$= \frac{-q^2\left(\frac{\partial f(x)}{\partial x_i} - v_i\right)}{p\left\{\beta p - \sqrt{\beta^2 p^2 + q^2\left(\frac{\partial f(x)}{\partial x_i} - v_i\right)^2}\right\}}$ (2.8)

Similarly

$\frac{\beta p - \sqrt{\beta^2 p^2 + q^2\left(\frac{\partial f(x)}{\partial x_i} - v_i\right)^2}}{p\left(\frac{\partial f(x)}{\partial x_i} - v_i\right)} = \frac{\beta^2 p^2 - \beta^2 p^2 - q^2\left(\frac{\partial f(x)}{\partial x_i} - v_i\right)^2}{p\left(\frac{\partial f(x)}{\partial x_i} - v_i\right)\left\{\beta p + \sqrt{\beta^2 p^2 + q^2\left(\frac{\partial f(x)}{\partial x_i} - v_i\right)^2}\right\}}$

$= \frac{-q^2\left(\frac{\partial f(x)}{\partial x_i} - v_i\right)}{\left\{\beta p + \sqrt{\beta^2 p^2 + q^2\left(\frac{\partial f(x)}{\partial x_i} - v_i\right)^2}\right\}}$ (2.9)

For $0 < x_i < \frac{q}{p}$ it is clear that

$$p^2 x_i^2 - q^2 = -(px_i + q)(q - px_i) < 0$$

and that for any point $x_i \left(\frac{\partial f(x)}{\partial x_i} - v_i\right)$ can only be zero, negative or positive.

**Case 1:**

Considering $\left(\frac{\partial f(x)}{\partial x_i} - v_i\right) = 0$ then $d_i(x, \lambda) = 0$

From $d_i(x, \lambda) - px_i < 0$, we obtain $0 < x_i$

Thus $(q - px_i) < (px_i + q)$

So, $-\left(\frac{p}{px_i + q} - \frac{p}{q - px_i}\right) > 0$

and $\frac{\partial}{\partial x_i}(L(x, \lambda)) = \frac{\partial f(x)}{\partial x_i} - v_i - \beta\left(\frac{p}{px_i + q} - \frac{p}{q - px_i}\right) > 0$

**Case II:**

Let $\left(\frac{\partial f(x)}{\partial x_i} - v_i\right) < 0$

From $d_i(x, \lambda) - x_i < 0$

Using (2.8) and (2.9) we have

$$0 < \left\{ x_i + \frac{\beta p + \sqrt{\beta^2 p^2 + q^2\left(\frac{\partial f(x)}{\partial x_i} - v_i\right)^2}}{p\left(\frac{\partial f(x)}{\partial x_i} - v_i\right)} \right\}$$

$$\left\{ x_i + \frac{\beta p - \sqrt{\beta^2 p^2 + q^2\left(\frac{\partial f(x)}{\partial x_i} - v_i\right)^2}}{p\left(\frac{\partial f(x)}{\partial x_i} - v_i\right)} \right\}$$

$$= \left\{ x_i - \frac{q^2\left(\frac{\partial f(x)}{\partial x_i} - v_i\right)}{p\left\{\beta p - \sqrt{\beta^2 p^2 + q^2\left(\frac{\partial f(x)}{\partial x_i} - v_i\right)^2}\right\}} \right\}$$

$$\left\{ x_i - \frac{q^2\left(\frac{\partial f(x)}{\partial x_i} - v_i\right)}{p\left\{\beta p + \sqrt{\beta^2 p^2 + q^2\left(\frac{\partial f(x)}{\partial x_i} - v_i\right)^2}\right\}} \right\}$$

$$x_i^2 - \frac{q^2 x_i\left(\frac{\partial f(x)}{\partial x_i} - v_i\right)}{p\left\{\beta p - \sqrt{\beta^2 p^2 + q^2\left(\frac{\partial f(x)}{\partial x_i} - v_i\right)^2}\right\}}$$

$$- \frac{q^2 x_i\left(\frac{\partial f(x)}{\partial x_i} - v_i\right)}{p\left\{\beta p + \sqrt{\beta^2 p^2 + q^2\left(\frac{\partial f(x)}{\partial x_i} - v_i\right)^2}\right\}}$$

$$+ \frac{q^4\left(\frac{\partial f(x)}{\partial x_i} - v_i\right)^2}{p^2\left\{\beta^2 p^2 - \beta^2 p^2 - q^2\left(\frac{\partial f(x)}{\partial x_i} - v_i\right)^2\right\}}$$

$$= x_i^2 - \frac{q^2 x_i\left(\frac{\partial f(x)}{\partial x_i} - v_i\right)}{p}\left\{\frac{1}{\beta p - \sqrt{\beta^2 p^2 + q^2\left(\frac{\partial f(x)}{\partial x_i} - v_i\right)^2}} + \frac{1}{\beta p + \sqrt{\beta^2 p^2 + q^2\left(\frac{\partial f(x)}{\partial x_i} - v_i\right)^2}}\right\} - \frac{q^2}{p^2}$$

$$= x_i^2 - \frac{q^2 x_i\left(\frac{\partial f(x)}{\partial x_i} - v_i\right)}{p}\left\{\frac{\beta p + \sqrt{\beta^2 p^2 + q^2\left(\frac{\partial f(x)}{\partial x_i} - v_i\right)^2} + \beta p - \sqrt{\beta^2 p^2 + q^2\left(\frac{\partial f(x)}{\partial x_i} - v_i\right)^2}}{\beta^2 p^2 - \beta^2 p^2 - q^2\left(\frac{\partial f(x)}{\partial x_i} - v_i\right)^2}\right\} - \frac{q^2}{p^2}$$

$$= x_i^2 - \left\{\frac{q^2 x_i\left(\frac{\partial f(x)}{\partial x_i} - v_i\right)}{p}\right\}\left\{\frac{2\beta p}{-q^2\left(\frac{\partial f(x)}{\partial x_i} - v_i\right)^2}\right\} - \frac{q^2}{p^2}$$

$$= x_i^2 - \frac{q^2}{p^2} + \frac{2\beta x_i}{\left(\frac{\partial f(x)}{\partial x_i} - v_i\right)} = \frac{p^2 x_i^2 - q^2}{p^2} + \frac{2\beta x_i}{\left(\frac{\partial f(x)}{\partial x_i} - v_i\right)}$$

$$\Rightarrow (p^2 x_i^2 - q^2)\left(\frac{\partial f(x)}{\partial x_i} - v_i\right) + 2\beta x_i p^2 < 0$$

$$\Rightarrow \frac{\partial f(x)}{\partial x_i} - v_i + \frac{2\beta p^2 x_i}{p^2 x_i^2 - q^2} > 0$$

because $\frac{\partial f(x)}{\partial x_i} - v_i < 0$

$$\Rightarrow \frac{\partial f(x)}{\partial x_i} - v_i - \beta\left[\frac{p}{px_i + q} - \frac{p}{q - px_i}\right] > 0$$

$$\Rightarrow \frac{\partial L(x, \lambda)}{\partial x_i} > 0$$

**Case III:**

If we have

$$0 > \left\{ x_i + \frac{\beta p + \sqrt{\beta^2 p^2 + q^2\left(\frac{\partial f(x)}{\partial x_i} - v_i\right)^2}}{p\left(\frac{\partial f(x)}{\partial x_i} - v_i\right)} \right\}$$

$$\left\{ x_i + \frac{\beta p - \sqrt{\beta^2 p^2 + q^2\left(\frac{\partial f(x)}{\partial x_i} - v_i\right)^2}}{p\left(\frac{\partial f(x)}{\partial x_i} - v_i\right)} \right\}$$

$$= x_i^2 - \frac{q^2}{p^2} + \frac{2\beta x_i}{\left(\frac{\partial f(x)}{\partial x_i} - v_i\right)} = \frac{p^2 x_i^2 - q^2}{p^2} + \frac{2\beta x_i}{\left(\frac{\partial f(x)}{\partial x_i} - v_i\right)}$$

Thus $0 > (p^2 x_i^2 - q^2)\left(\frac{\partial f(x)}{\partial x_i} - v_i\right) + 2\beta x_i p^2$

Therefore $0 < \frac{\partial f(x)}{\partial x_i} - v_i + \frac{2\beta p^2 x_i}{p^2 x_i^2 - q^2}$

$-\frac{\partial f(x)}{\partial x_i} - v_i - \beta \left[\frac{2}{px_i + q} - \frac{p}{q - px_i}\right] - \frac{\partial L(x,\lambda)}{\partial x_i}$

Using this result the rest of this Lemma can be proved similarly.

Now if we will consider all $-\frac{q}{p} < x < \frac{q}{p}$

$d(x, \lambda) - x = 0$ if and only if $\nabla_x L(x, \lambda) = 0$

Since $-\frac{q}{p} < d(x, \lambda) < \frac{q}{p}$ then $d(x, \lambda) - x$ satisfies the desired property and when searching for a point in $(x, \lambda) - x$, the constraint $-\frac{q}{p} < x < \frac{q}{p}$ is always satisfied automatically if the step length is a number between $-\frac{q}{p}$ and $\frac{q}{p}$.

Let $F(\lambda) = -\left(\frac{q\left(\frac{\partial f(x)}{\partial x_i} - \lambda\right)}{\beta p + \sqrt{\beta^2 p^2 + q^2\left(\frac{\partial f(x)}{\partial x_i} - \lambda\right)^2}} + \frac{q\left(\frac{\partial f(x)}{\partial x_i} - \lambda\right)}{\beta p + \sqrt{\beta^2 p^2 + q^2\left(\frac{\partial f(x)}{\partial x_i} - \lambda\right)^2}}\right)$

Given any point $x \in F$ with $-\frac{q}{p} < x < \frac{q}{p}$ for $d(x, \lambda) - x$ to become a feasible descent direction as $h(x, \beta) = f(x) + \beta b(x)$ and we want to solve $p(\lambda) = 0$ (2.10)

$\lambda(x) = \frac{1}{2q}\left(\frac{\partial f(x)}{\partial x_s} + \frac{\partial f(x)}{\partial x_t}\right)$ is a solution of (2.10)

Since $\{\nabla f(x) / x \in \beta\}$ is bounded hence $\{\lambda(x) : x \in B\}$ is bounded. Using the feasible descent direction $d(x, \lambda(x)) - x$ and $\lambda(x) = \frac{1}{2q}\left(\frac{\partial f(x)}{\partial x_s} + \frac{\partial f(x)}{\partial x_t}\right)$ for updating

Lagrange multiplier λ we have developed an algorithm for approximating a solution of

$\min : f(x) - \frac{1}{2} x^T (W - \alpha I) x$

subject to $x_s + x_t = 0$

$x_i \in \{-1, 1\}, i = 1, 2, \ldots, n$

Let $\beta_m, m = 1, 2, \ldots$ be any given sequence of positive numbers such that $\beta_1 > \beta_2 > \cdots$ and $\lim_{m \to \infty} \beta_m = 0$. The value of $\beta_1$ should be sufficiently large so that $h(x, \beta_1)$ is convex over $-\frac{q}{p} < x < \frac{q}{p}$.

Let $x^0$ be an arbitrary non zero interior point of B satisfying $e_{it}^T x^0 = 0$. For m=1, 2,… starting at $x^{m-1}$, we can generate the feasible descent direction $d(x, \lambda(x)) - x$ to search for better feasible interior point $x^m$ satisfying $d(x^m, \lambda(x^m)) = 0$

## 3 Main Results

Using the above lemmas we can prove the following main theorem.

**Theorem 1.**

For $\beta = \beta_k$ every limit point of $x^k, k = 1, 2, \ldots$, generated by $x^{k+1} = x^k + \mu_k(d(x^k, \lambda(x^k)) - x^k)$ is a stationary point of $\min h(x, \beta) = f(x) + \beta b(x)$ subject to $x_s + x_t = 0$

**Proof:** $\beta_k, k = 1, 2, \ldots$ be any given sequence of positive numbers such that $\beta_1 > \beta_2 > \cdots$ and $\lim_{k \to \infty} \beta_k = 0$. The value of $\beta_1$ should sufficiently large so that $h(x, \beta)$ is convex over $-\frac{q}{p} < x < \frac{q}{p}$.

Let $x^0$ be an arbitrary non zero interior point of B satisfying $e_{it}^T x^0 = 0$ for k=1,2,… starting from $x^{k-1}$, we derive the feasible descent direction $d(x^k, \lambda(x^k)) - x^k = 0$, if $\|d(x^k, \lambda(x^k)) - x^k\| < 1$

When β is sufficiently small so that a feasible solution which every component equal to either -1 or 1 can be generated by rounding off $x^k$.

Let $\beta_{k+1} = \theta \beta_k$

$x^{k+1} = x^k + \mu_k(d(x^k, \lambda(x^k)) - x^k))$

Satisfying

$h(x^{k+1}, \beta_k) = \min_{\mu \in [0,1]} h[x^k + \mu(d(x^k, \lambda(x, k)) - x^k, \beta_k)]$

Let $a_{max} = \max_{1 \leq i \leq n} \max_{x \in B} \left( \frac{\partial f(x)}{\partial x_i} - v_i(x) \right)$

with

$x^{min} = (x_1^{min}, x_2^{min}, \ldots, x_n^{min})^T$,
$x^{max} = (x_1^{max}, x_2^{max}, \ldots, x_n^{max})^T$, $\mu[0,1]$

$x_i^{min} = \min \left\{ x_i^0, \frac{-a_{max}}{\beta p + \sqrt{\beta^2 p^2 + a_{max}^2}} \right\}$

$x_i^{max} = \max \left\{ x_i^0, \frac{-a_{max}}{\beta p + \sqrt{\beta^2 p^2 + a_{max}^2}} \right\}, i = 1,2,\ldots n$

Since $\{\nabla f(x): x \in B\}$ and $\{\lambda(x): x \in B\}$ are bounded then $a_{max}$ is bounded.

Thus $-\frac{p}{q} < x^{min} \leq d(x, \lambda(x)) \leq x^{max} < \frac{p}{q}$

For any $x \in B$ and $x^{min} \leq x^k \leq x^{max}, k=1,2,\ldots$

Therefore no limit of $x_i^{\beta}, k = 1,2\ldots$ is equal to either -1 or 1 for i=1,2,…n

From the above lemma we obtain $d(x^k, \lambda(x^k)) - x^k$ is a feasible solution of $\min h(x, \beta) = f(x) + \beta b(x)$ subject to $x_s + x_t = 0$

Let $X = \{x \in F: x^{min} \leq x^k \leq x^{max}\}$

and $\Omega = \{x \in X: d(x, \lambda(x)) - x = 0\}$, for any $x \in X$

Let $A(x) = \{x + \mu^*(d(x, \lambda(x)) - x)\}$

$\mu^* \in [0,1], h(x + \mu^*(d(x, \lambda(x)) - x), \beta)$

$= \min_{\mu \in [0,1]} h(x + \mu(d(x, \lambda(x)) - x), \beta)$

Now we will prove A(x) is closed at every point $x \in X - \Omega$

Let $\bar{x}$ be any arbitrary point of $X - \Omega$

Let $x^r \in X - \Omega, r = 1,2,\ldots$ be a sequence convergent to $\bar{x}$ and $y^r \in A(\bar{x})$, r=1,2,…n be a convergent sequence converges to $\bar{y}$.

To prove that $A(\bar{x})$ is closed, we need to show that $\bar{y} \in A(\bar{x})$. From $x^r \in X - \Omega$ and $\bar{x} \in X - \Omega$ we have $d(x^r, \lambda(x^r)) - x^r \neq 0$ and $(\bar{x}, \lambda(\bar{x})) - \bar{x} \neq 0$, where $d(x, \lambda(x))$ is continuous. Thus $d(x^r, \lambda(x^r))$ converges to $d(\bar{x}, \lambda(\bar{x}))$ as $r \to \infty$. Since $y^r \in A(x^r)$ then there exists a number $\mu^* \in [0,1]$ satisfying

$y^r = x^r + \mu^* d(x^r, \lambda(x^r)) - x^r \in A(x^r)$

From $(x^r, \lambda(x^r)) - x^r \neq 0$, we obtain

$\mu_r^* = \frac{\|y^r - x^r\|}{d\|x^r, \lambda(x^r) - x^r\|}$

And as, $r \to \infty$

$\mu^* - \overline{\mu^*} = \frac{\|\bar{y} - \bar{x}\|}{d\|\bar{x}, \lambda(\bar{x})) - \bar{x}\|}$

with $\bar{\mu} \in [0,1], \bar{y} = \bar{x} + \overline{\mu^*}[d(\bar{x})(\lambda(\bar{x})) - \bar{x}]$, since $y^r \in A(x^r)$

we have

$h(y^r, \beta) \leq h x^r + \mu(d(x^r, \lambda(x^r)) - x^r), \beta)$ for any $\mu \in [0,1]$

It implies that $h(\bar{y}, \beta) \leq h(\bar{x} + \mu(d(\bar{x}) \cdot (\lambda(\bar{x})) \quad \bar{x})), \beta)$

For any $\mu \in [0,1]$, which proves that h

$(\bar{y}, \beta) = \min h(\bar{x} + \mu(d(\bar{x}) \cdot \lambda(\bar{x})) - \bar{x}), \beta)$

$\mu \in [0,1]$

According to the definition of A (x), it follows that $\bar{y} \in A(\bar{x})$, since X is bounded and $x^k \in X, k = 1,2,\ldots$

Then by Bolzonoweierstrass theorem we can extract a convergent subsequence from the sequence

$x^k, k = 1,2,\ldots$

Let $x^{k_i}, i = 1,2,\ldots$ be a convergent subsequence of the sequence $x^k, k = 1,2,\ldots$

Let $x^*$ be the limit point of the sequence. Clearly as $k \to \infty$, $h(x^k, \beta)$ converges to $h(x^*, \beta)$, since $h(x, \beta)$ is continuous and $h(x^{k+1}, \beta) < h(x^k, \beta), k = 1,2,\ldots$

Considering the sequence $x_j^{k+1}, j = 1,2,\ldots$ and by (2.10)

We can write $x_j^{k+1} = x^{k_i} + \mu_{k,j}(d(x^{k_i}, \lambda(x^{k_i})) - x^{k_i})$

and $h(x_j^{k+1}, \beta) = \min h(x^{k_i} + \mu(d(x^{k_i}, \lambda)(x^{k_i} - x^{k_i}), \beta)$

According to definition of $A(x)$ we have $x_j^{k+1} \in A(x^{k_i}), j = 1,2,\ldots$ are bounded, then there exist $x_j^{k+1}, j \in k$ a convergent subsequence. If $x^* \in \Omega$ and $A(x^*)$ is closed then $x^* \in A(x^*)$ and we have $h(x^*, \beta) < h(x^k, \beta)$

Which contradicts that $h(x^k, \beta)$ converges as $\to \infty$ i.e $x^* \in \Omega$.

The use of logarithmic barrier function finds a minimal point for solving linearity constrained continuous optimization problem to s-t max cut problem.

## 4 Numerical Example

In order to establish the effectiveness and efficiency of the algorithm for obtaining optimization, we have solved by using MATLAB.

For our solution we have consider $\alpha$ = 0.000001

$\beta_0$ = 1- $S_{min}/2$ with $S_{min}$ being the minimum eigen value of W-$\alpha$I.

In our computation we have taken following variables.

NI=No. of iteration

OBJM =Objective value of (1.4)

OBJU =Greatest integer value of

$$\max : \frac{1}{4}\sum_{ij} w_{ij}(1 - x_{ij})$$

Subject to trace $(e_{it}e_{it}^T X) = 0$

$X_{ii} = 1$, $i = 1,2,...n$

$x \geq 0$, RATIO= (OBJU-OBJM)/OBJU

In the computation we have considered weights $w_{ij}$ are the random integer between 1 to 50 and the results computed are given in the following table for the value of m=.6

Table for the numerical results for m=.6

| No of Test | NI  | OBJM  | OBJU  | RATIO |
|------------|-----|-------|-------|-------|
| 1          | 100 | 35670 | 35798 | .03   |
| 2          | 150 | 36150 | 37630 | .02   |
| 3          | 200 | 46260 | 47350 | .02   |
| 4          | 250 | 20350 | 21450 | .02   |
| 5          | 300 | 35672 | 36505 | .02   |

From the above computation it shows that the ratio is nearly equal to .02 and it proves the convergence of the solution to the local minimum point and shows that the algorithm is an effective one.

## 5 Conclusion

In this paper we have taken a general logarithmic barrier function to solve the continuous optimization problem by transforming it in to s-t max-cut problem and developed an algorithm is based on the barrier generalized logarithm the barrier function. The algorithm developed for generating decreasing sequence of solution which converges at least a local minimum point of (2.1) with every component equal to either $\frac{-y}{p}$ $to$ $\frac{y}{p}$.

**Dr.A.K.Ojha:** Dr A.K.Ojha received a Ph.D (Mathematics) from Utkal University in 1997.Currently he is an Asst. Prof. in Mathematics at I.I.T Bhubaneswar, India. He is performing research in Neural Network, Genetical Algorithm, Geometric Programming and Particle Swarm Optimization. He is served more than more than 27 years in different Govt. colleges in the state of Orissa. He is published 22 research paper in different journals and 7 books for degree students such as: Fortran 77 Programming, A text book of Modern Algebra, Fundamentals of Numerical Analysis etc.

**Dr.C.Mallick:** Dr.C.Mallick received a Ph.D (Mathematics) from Utkal University in 2008.Currently he is an Lecturer in Mathematics at B.OS.E Cuttack, India. He is performing research in Neural Network. He is published 3 books for degree students such as: A Text Book of Engineering Mathematics, Interactive Engineering Mathematics etc.

**D. Mallick:** Mr.D.Mallick received a M.Phil (Mathematics) from Sambalpur University in 2002.Currently he is an Asst. Prof. in Mathematics at Centurion Institute of Technology, Bhubaneswar, India. He is performing research in Neural Network and Optimization Theory. He is served more than more than 6 years in different Colleges in the state of Orissa. He is published 6 research paper in different journals.